\documentstyle[12pt]{article}

\begin{document}
\newcommand{\ol }{\overline}
\newcommand{\ul }{\underline }
\newcommand{\ra }{\rightarrow }
\newcommand{\lra }{\longrightarrow }
\newcommand{\ga }{\gamma }
\newcommand{\st }{\stackrel }
\newcommand{\scr }{\scriptsize }
\title{\Large\textbf{  Some Properties of the Schur Multiplier
with Algebraic Topological Approach }}
\author{\textbf{
  Behrooz Mashayekhy
\footnote{Correspondence: mashaf@math.um.ac.ir} \ \ and Hanieh
Mirebrahimi } \\
Department of Mathematics,\\ Center of Excellence in Analysis on Algebraic Structures,
\\ Ferdowsi University of Mashhad,\\
P. O. Box 1159-91775, Mashhad, Iran}
\date{ }

\maketitle

\begin{abstract}

In this paper, using a relationship between the Schur multiplier of
a group $G$, the fundamental group, and the second homology group of
the Eilenberg-MacLane space of $G$, we present new proofs for some
famous properties of the Schur multiplier and its structure for the
free product and the direct product with an algebraic topological
approach. Also, we try to find the Schur multiplier of a free
amalgamated product in terms of Schur multipliers of its factors,
which is a new result.
\end{abstract}
\ \ \\
A. M. S. classification: 57M07; 55Q20; 55P20; 20C25.
\newpage
\hspace{-0.65cm}\textbf{1. Introduction and Preliminaries}\\

Let $F/R$ be a free presentation of a group $G$, then the Schur
multiplier of $G$ is defined to be $$M(G)=(R\cap F')/[R,F]$$ ( see
[2] for further details ). There are some well-known facts about
$M(G)$. By a theorem of Schur [2], with respect to the
presentation of $G$, there exists an upper bound for the number of
the generators of $M(G)$. In particular, if $G$ is a cyclic group,
$M(G)$ is trivial [2]. A theorem of C. Miller [2], says that
$$M(G_1*G_2)\cong M(G_1)\oplus M(G_2).$$ Also Schur Theorem
[2] asserts that $$M(G_1\times G_2)\cong M(G_1)\oplus M(G_2)\oplus
G_1\otimes G_2.$$

In this paper, we are going to give some new proofs for the above
facts with algebraic topological methods. Also, as a new result,
we try to find the Schur multiplier of a free amalgamated product
in terms of Schur multipliers of its factors. We suppose that the
reader is familiar with some well-known notions in algebraic
topology such as homotopy groups, homology groups, CW-spaces, and
basic notions in group theory.

In order to find a suitable relationship between the Schur
multiplier and some famous notions of algebraic topology, we
mention the following notes.

\hspace{-0.65cm}\textbf{\bf Theorem 1.1 (Hurewicz Theorem [4]).}
For any path connected space $X$, the Hurewicz map induces an
isomorphism called Hurewicz isomorphism, as follows:
$$\frac{\pi_1(X,x_0)}{(\pi_1(X,x_0))'}\cong H_1(G).$$

\hspace{-0.65cm}\textbf{\bf Theorem 1.2 (Hopf Formula [1]).} If
$K$ is a CW-complex with $\pi_1(K)=G$ and $F/R$ is a presentation
of $G$ then there exists the following isomorphism along the
Hurewicz map:
$$\frac{H_2(K)}{h_2(\pi_2(K))}\cong \frac{R\cap F'}{[R,F]},$$
where $H_2(K)$ is the second homology group of $K$ and $h_2$ is
the corresponding Hurewicz map.

\hspace{-0.65cm}\textbf{\bf Theorem 1.3 ([6]).} For any group
$G$, there exists a CW-complex $X$ with $$\pi_1(X)\cong G\ and\
\pi_n(X)=1\ for\ all\ n\geq 2.$$ The space $X$ is called
Eilenberg-MacLane space of $G$.

\hspace{-0.65cm}\textbf{\bf Corollary 1.4.} For any group $G$ and
its Eilenberg-MacLane space, $K$ say, we have
$$\pi_1(K)\cong G\ and\ H_2(K)\cong M(G). $$

\hspace{-0.65cm}\textbf{\bf Remark 1.5 ([4,6]).} For any
CW-complex $K$ with $k_i$ many $i$-cells, $d(H_i(K))\leq k_i$. In
addition, the number of $2$-cells in an Eilenberg-MacLane space,
obtained from a presentation of a group $G$, is equal to the
number of it's relators. (See the proof of [6, Note 6.44])

\hspace{-0.65cm}\textbf{\bf Theorem 1.6 ([3,5]).} For any
 numbers $n,m\in N$, there exists a CW-complex $L(n,m)$,
called Lens space, with $$\pi_1(L(n,m))\cong Z_n\ and\
H_2(L(n,m))=1.$$

\hspace{-0.65cm}\textbf{\bf Theorem 1.7 (Mayer-Vietoris Sequence
[4]).} For any two subcomplexes $X_1$ and $X_2$ of a CW-complex
$X$, with $X=X_1\cup X_2$, there is an exact sequence as follows
$$\cdots\rightarrow H_n(X_1\cap X_2)\rightarrow H_n(X_1)\oplus
H_n(X_2)\rightarrow H_n(X)\rightarrow H_{n-1}(X_1\cap
X_2)\rightarrow\cdots.$$

\hspace{-0.65cm}\textbf{\bf Theorem 1.8 (Kunneth Formula [4]).}
For any pair of spaces $X$ and $Y$ and for every integer $n\geq
0$, we have the following relation between their homology groups
$$H_n(X\times Y)\cong \sum_{i+j=n}H_i(X)\otimes H_i(Y)\oplus
\sum_{p+q=n-1} Tor(H_p(X),H_q(Y)).$$

\hspace{-0.65cm}\textbf{2. Main Results }\\

In this section, first, using the previous relationship between
the Schur multiplier, homology and homotopy groups, we try to give
algebraic topological proofs for some well-known facts about the
Schur multiplier.\\

\hspace{-0.65cm}\textbf{\bf Theorem 2.1.} \textit{For any
k-relator group $G$, there exists an upper bound for the number of
the generators of its Schur multiplier, as follows}
$$d(M(G))\leq k.$$

\hspace{-0.7cm}\textbf{ Proof.} First, suppose that $K$ is the
Eilenberg-MacLane space of $G$, then using Remark 1.5, the number
of $2$-cells in the complex $K$ equals exactly to $k$ and
consequently $d(H_2(K))\leq k$. Finally by Hopf Formula for the
Eilenberg-MacLane space $K$, $ M(G)\cong H_2(K)$ and so the result
holds. $\Box$\\

\hspace{-0.65cm}\textbf{\bf Theorem 2.2.} \textit{The Schur
multiplier of any cyclic group $G$ is trivial.}\\

\hspace{-0.7cm}\textbf{ Proof.} If $G$ is an infinite cyclic
group, we can consider a circle $S^1$ as it's Eilenberg-MacLane
space ($S^1$ is a CW-complex whose fundamental group is infinite
cyclic and it's higher homotopy groups are trivial [4]). Hence
using the Hopf isomorphism $M(G)\cong H_2(S^1)$ and the fact of
$S^1$ whose second homology group is trivial, in this case, the
result is satisfied.

 Also for a finite cyclic group $G$ of order
$n$, using Theorem 1.2, we have the Lens space $L(n,1)$ as an
Eilenberg-MacLane space of $G$. So by a similar argument to the
above and Theorem 1.2, the proof is completed. $\Box$\\

\hspace{-0.65cm}\textbf{\bf Theorem 2.3.} \textit{For any two
groups $G_1$ and $G_2$, we have the isomorphism}
$$M(G_1*G_2)\cong M(G_1)\oplus M(G_1).$$

\hspace{-0.7cm}\textbf{ Proof.} First, using Theorem 1.1, let
$K_1$ and $K_2$ be the Eilenberg-MacLane spaces of $G_1$ and
$G_2$, respectively. By Van-Kampen Theorem for the fundamental
group of wedge space, $\pi_1(K_1\vee K_2)\cong
\pi_1(K_1)*\pi_1(K_2)$. Also using definitions of $\pi_n$ and
wedge space, $\pi_n(K_1\vee K_2)=1$, for all $n\geq 2$. Hence the
wedge spaces $K_1\vee K_2$ can be considered as an
Eilenberg-MacLane spaces of $G_1*G_2$ and with respect to the Hopf
Theorem, we have
$$ M(G_1)\cong H_2(K_1)\ ,\ M(G_2)\cong H_2(K_2)$$$$\&\ \ M(G_1*G_2)\cong H_2(K_1\vee K_2).$$

Finally, using the Mayer-Vietories sequence for the wedge space
$K_1\vee K_2$, we conclude the following isomorphism which
completes the proof
$$ H_2(K_1\vee K_2)\cong H_2(K_1)\oplus
H_2(K_2).\ \Box$$

\hspace{-0.65cm}\textbf{\bf Theorem 2.4.} \textit{For any two
groups $G_1$ and $G_2$, the following isomorphism holds:}
$$M(G_1\times G_2)\cong M(G_1)\oplus M(G_2)\oplus
(G_1)_{ab}\otimes (G_2)_{ab}.$$

\hspace{-0.7cm}\textbf{ Proof.} Similar to the previous proof,
suppose that $K_1$ and $K_2$ are the Eilenberg-MacLane spaces of
$G_1$ and $G_2$, respectively. Using one of the properties of the
homotopy functor $\pi_n$ to preserve the direct product, we have
$\pi_1(K_1\times K_2)\cong \pi_1(K_1)\times\pi_1(K_2)$ and
$\pi_n(K_1\times K_2)=1$, for all $n\geq 1$.

So the space $K_1\times K_2$ is an Eilenberg-MacLane spaces of
$G_1\times G_2$. Hence by  Hopf isomorphism, we have $M(G_i)\cong
H_2(K_i)$, for $i=1,2$, and
$$M(G_1\times G_2)\cong H_2(K_1\times K_2).$$
Also, using Kunneth Formula, the properties of the functor Tor and
tensor product, and the fact of $H_0(X)$ which is isomorphic to
the infinite cyclic group $Z$ where $X$ is a path connected space
( note that Eilenberg-MacLane spaces are path connected ), we
have the following relation between the first and the second
homology groups,
$$H_2(K_1\times K_2)\cong H_2(K_1)\oplus H_2(K_2)\oplus H_1(K_1)\otimes
H_1(K_2).$$ Finally, by Hurewicz isomorphisms $H_1(K_i)\cong
(\pi_1(K_i))_{ab}=(G_i)_{ab}$, for $i=1,2$, we conclude the
result of the theorem. $\Box$\\

Now, as a new result, we try to find the structure of the Schur
multiplier of a free amalgamated product in terms of the Schur
multipliers of its factors.\\

\hspace{-0.65cm}\textbf{\bf Theorem 2.5.} \textit{Suppose that
$G$ is the free amalgamated product of its two subgroups $G_1$
and $G_1$ over a subgroup $H$. Then we have the following exact
sequence, which is induced by the Mayer-Vietories sequence along
the Hurewicz map:}
$$\cdots \rightarrow M(H)\rightarrow M(G_1)\oplus M(G_2)\rightarrow M(G)
$$ $$\rightarrow H_{ab} \rightarrow G_{1 ab}\oplus G_{2
ab}\rightarrow G_{ab}\rightarrow \cdots .$$

\hspace{-0.7cm}\textbf{ Proof.} First, we consider an
Eilenberg-Maclane space corresponding to the presentation of $H$,
$Y$ say, and note that we can extend the presentation of $H$ to a
presentation for $G_1$ and a presentation for $G_2$.

Also, by joining some $1$-cells and attaching $2$-cells via the
relations, similar to the method of [4, Theorem 7.45] and [6,
Note 6.44], we extend $Y$ to Eilenberg-Maclane spaces $X_1$ and
$X_2$ corresponding to the above presentations of $G_1$ and $G_2$,
respectively. Note that the construction is considered in such a
way that $Y$ is a deformation retract of the space $X_1\cap X_2$.

Now using the van-Kampen theorem, the fundamental group
$\pi_1(X_1\cup X_2)$ is the free amalgamated product of two groups
$\pi_1(X_1)\cong G_1$ and $\pi_1(X_2)\cong G_2$ over the subgroup
$\pi_1(Y)\cong \pi_1(X_1\cap X_2)\cong H$ [5].

Hence by uniqueness of the free amalgamated product up to
isomorphism, we conclude that
$$G\cong \pi_1(X_1\cup X_2).$$

Also, we have the Mayer-Vietoris Sequence for CW-complexes $X_1$
and $X_2$ [4], as follow:

$$\cdots\rightarrow H_2(X_1\cap X_2)\rightarrow H_2(X_1)\oplus
H_2(X_2)\rightarrow H_2(X_1\cup X_2)$$ $$\rightarrow H_1(X_1\cap
X_2)\rightarrow H_1(X_1)\oplus H_1(X_2)\rightarrow \rightarrow
H_2(X_1\cup X_2)\rightarrow 1.$$

Finally, using Hurewicz maps and so Hopf Formula (See the proof of
[1, Theorem(Hopf)]), we get the
result. $\Box$\\

As some corollaries of the recent theorem, we are going to present
the following notes about the structure of the Schur multiplier of
a free amalgamated product in some special cases.\\

\hspace{-0.65cm}\textbf{\bf Remark 2.6.} Let $G$ be the free
amalgamated product of its two subgroups $G_1$ and $G_2$ over a
subgroup $H$. Then by the exact sequence in previous theorem, the
following hold.

i) If the Schur multiplier of $H$ is trivial, then $M(G)$ contains
a subgroup isomorphic to the group $M(G_1)\oplus M(G_2)$.

ii) If $H$ is a perfect group with trivial Schur multiplier, then
$M(G)$ is isomorphic to the group $M(G_1)\oplus M(G_2)$.

\end{document}